\documentclass[12pt]{amsart}
\usepackage[francais]{babel}
\usepackage{pdfpages}
\usepackage[T1]{fontenc}
\usepackage{graphics}
\usepackage{graphicx}

\title{A propos des nombres de Conway :
\\ Lettre \`a un ami} 

\author{Labib Haddad}
\address{120 rue de Charonne, 75011 Paris, France}
\email{labib.haddad@wanadoo.fr}

\usepackage{amssymb}
\usepackage{amsmath}

\newcommand{\su}{\subsection*}
\newcommand{\head}{\section*}
\newcommand{\noi}{\noindent}

\newcommand{\Ž}{\'e}

\newcommand{\ˆ}{\`a}
\newcommand{\}{\`u}

\newcommand{\D}{\mathbb D}

\newcommand{\N}{\mathbb N}
\newcommand{\Q}{\mathbb Q}
\newcommand{\R}{\mathbb R}

\newcommand{\Z}{\mathbb Z}

\newcommand{\cal}{\mathcal}

\newcommand{\sA}{\mathsf A}
\newcommand{\sB}{\mathsf B}
\newcommand{\sO}{\mathsf O}
\newcommand{\sM}{\mathsf M}
\newcommand{\sN}{\mathsf N}

\newcommand{\sS}{\mathsf S}

\newcommand{\sX}{\mathsf X}
\newcommand{\sY}{\mathsf Y}

\newcommand{\f}{\varphi}

\newcommand{\leqs}{\leqslant}

\newcommand{\geqs}{\geqslant}

\newcommand{\guil}{\guillemotleft}  
\newcommand {\guir}{\guillemotright}

\newcommand {\et}{\ \text{et}\ }

\newcommand {\ou}{\ \text{ou}\ }

\newcommand {\si}{\ \text{si}\ }

\newcommand {\pourtout}{\ \text{pour tout}\ }

\newcommand {\pourtous}{\ \text{pour tous}\ }

\newcommand{\stm}{\smallsetminus}

\newcommand{\vide}{\emptyset}

\newcommand{\inc}{\subset}

\newcommand{\lopar}{\noi \{$\looparrowright$ \ }
\begin{document}
\maketitle
\thispagestyle{empty}

\markboth{Labib Haddad}{Nombres de Conway}

\centerline{\it \`A la m\Žmoire de John Horton Conway, 1937-2020}

\

\

Cher Ami,

\

Tu m'\Žcris ceci. 

{\it C'est la premi\re fois que j'entends parler des nombres de Conway. Pourrais-tu m'expliquer de quoi il s'agit, s'il te pla\"t ?
J'ai lu deux ou trois textes \ˆ ce sujet mais je n'y ai pas compris grand chose. Je crains d'avoir perdu un peu le contact avec les math\Žmatiques de notre temps!}

Rassure-toi, il n'en est rien !

Sache qu'il existe plusieurs mani\res diff\Žrentes, mais \Žquivalentes, pour introduire les nombres de Conway. Elles sont plus ou moins claires, plus ou moins adroites. Certaines peuvent \tre un peu absconses et (par coquetterie) deviennent lapidaires, aux d\Žpends de la clart\Ž. Pour ma part, je pense que la redondance est n\Žcessaire \ˆ la compr\Žhension. La r\Žp\Žtition est un bon outil p\Ždagogique. La concision est un luxe que l'on peut se payer seulement lorsque l'on a d\Žj\ˆ compris ...

\

\lopar Il y a deux r\Žcits de la cr\Žation dans la Gen\se et quatre \Žvangiles dans la Bible !\}

\

La mani\re de {\sc Conway} lui-m\me m'a sembl\Ž assez difficile \ˆ suivre. Elle se veut concise et ramass\Že  ! Cela n'en facilite pas la compr\Žhension. Cependant, certains la vante jusqu'au dithyrambe !

\

Une digression avant d'en venir \ˆ l'essentiel.

\

{\sc Henri Cartan}, tu t'en souviens, nous apprenait, directement ou par l'exemple, ceci. Lorsque l'on a d\Žmontr\Ž un \Žnonc\Ž donn\Ž et v\Žrifi\Ž son exactitude, il reste deux choses \ˆ faire. Se demander si l'on n'a pas, en fait, d\Žmontr\Ž davantage que l'\Žnonc\Ž lui-m\me (ce qui arrive parfois) puis voir si l'on ne peut pas simplifier cette d\Žmonstration et, ce faisant,  l'\Žtendre, plus avant.  Quant \ˆ {\sc Gustave Choquet}, il nous expliquait qu'en mati\re de g\Žn\Žralistion, il faut  savoir  garder la bonne distance, la bonne mesure, ni trop, ni trop peu !  Partant d'un cas pr\Žcis, si l'on g\Žn\Žralise trop, on se perd dans des banalit\Žs, si pas assez, on perd du recul !

Je vais essayer de te faire une pr\Žsentation des nombres de Conway, {\it \ˆ ma mani\re},  qui (\ˆ mon avis) satisfait les exigences de rigueur et de simplicit\Ž. 
\

\head{1 \ Un prologue}

\

\su{Rappel} Pour \Žviter les paradoxes de la th\Žorie naissante des ensembles \Žrig\Že par {\sc Cantor}, on a us\Ž de diverses axiomatisations. Tu connais celle de {\sc Zermelo} et {\sc Fraenkel}. Il  en est une, d\Žcrite en appendice dans le livre de Topologie g\Žn\Žrale de {\sc Kelley} o\ il pr\Žcise que c'est une variante  des syst\mes de {\sc Skolem} et de {\sc A.P. Morse}. Au c\™t\Ž des ensembles proprement dit, il y a les {\bf classes} que l'on d\Žsigne ainsi sous forme de classificateurs, $\{x : \f(x)\}$, ce qui se  lit comme suit : la classe des objets $x$ qui v\Žrifient la condition $\f(x)$.  Un ensemble est une classe, mais il y a des classes qui ne sont pas des ensembles ! Je ne vais pas m'\Žtendre davantage l\ˆ-dessus, tu as saisis, j'en suis s\žr. 

Chez {\sc Bourbaki}, un ensemble est d\Žfini par une relation collectivisante. Si la  relation n'est pas collectivisante elle d\Žfinit une classe qui n'est pas un ensemble. Dans la th\Žorie des ensembles, {\it \ˆ la Bourbaki},   les {\bf nombres ordinaux}, $0,1,2,\dots,\omega, \omega +1, \dots, \omega2, \omega2 + 1, \dots,$ ayant \Žt\Ž d\Žfini, on montre qu'ils ne forment pas un ensemble car cela introduirait une contradiction : on parle ainsi de {\bf la classe} des nombres ordinaux et pas de leur ensemble. Bien entendu, les ordinaux inf\Žrieurs \ˆ un ordinal donn\Ž $\alpha$ si grand soit-il forment bien un ensemble; l'ensemble $[0,\alpha[$ est un segment initial de la classe des ordinaux !

\

\noi Il en va de m\me des nombres de Conway : il s'agit d'une {\bf classe}, et pas d'un ensemble, une classe munie d'une structure de {\bf corps ordonn\Ž}, le plus grand possible, en un certain sens ! On d\Žfinit ces nombres par r\Žcurrence. A chaque \Žtape de la r\Žcurrence, on a bien {\bf un ensemble} de nombres de Conway, mais ils ne forment qu'un segment de la totalit\Ž !

\

Dans ma pr\Žsentation,  en place de l'expression \guil nombres de Conway\guir, je dirai plus simplement, et bri\vement, {\bf les nombres}. Je parlerai ainsi {\bf des nombres puis de leurs noms}.

\

\noi On aura besoin des pr\Žliminaires suivants.

\su{Coupures} Les coupures de Dedekind sont un outil bien connu qui  sert \ˆ construire le corps $\R$ des nombres r\Žels \ˆ partir du corps $\Q$ des nombres rationnels. On les trouvent dans la plupart des manuels d'analyse.

\

\noi En voici une g\Žn\Žralisation au cas des ensembles totalement ordonn\Žs quelconques.

\

\noi Soit $E$ un ensemble muni d'une relation $\leqs$ d'ordre total. Cela veut dire que $\leqs$ est une relation binaire sur l'ensemble $E$ : une relation qui  est r\Žflexive, transitive, antisym\Žtrique et totale. Autrement dit, quels que soient les \Žl\Žments $x,y,z,$ de $E$, on a les propri\Žt\Žs suivantes :

\

$x\leqs x$

$x\leqs y \et y\leqs z \implies x\leqs z$

$x \leqs y \et y\leqs x \implies x = y$

on a toujours $x\leqs y \ou y\leqs x$, l'un ou l'autre  ou les deux \ˆ la fois. 

[Il s'agit du \emph{ou} inclusif pas du \emph{ou} exclusif.]

\

\noi On \Žcrit $x < y$ lorsque l'on a ($x\leqs y \et x\neq y$). On dit que $<$ est {\bf l'ordre strict} associ\Ž \ˆ l'ordre $\leqs$.

\

\noi Bien entendu, $y\geqs x$ est synonyme de $x\leqs y$. De m\me, $y > x$ est synonyme de $x < y$.

\

\

\noi On appelle {\bf coupure} dans $E$, tout couple $(A | B)$ o\ $A$ et $B$ sont des parties de l'ensemble $E$ telles que, pour tous $a\in A$ et $b\in B$, on ait $a \leqs b$ et telles que l'on ait
$$A\cap B = \vide \ , \  A\cup B = E.$$

\su{Ordre total sur l'ensemble des coupures}  L'ensemble $\cal C(E)$ des coupures dans $E$ est muni, d'une mani\re naturelle, d'une relation d'ordre totale comme suit.  Etant donn\Žes  $(A|B)$ et $(C|D)$, deux coupures dans $E$, on \Žcrit $(A|B) \leqs (C|D)$ lorsque l'on a $A\inc C$. On v\Žrifie [il faut le faire soi-m\me] que c'est bien une relation d'ordre totale sur l'ensemble $\cal C(E)$.

\

\su{Amalgame} Soit  alors $\cal A = E\cup \cal C(E)$, la r\Žunion de l'ensemble totalement ordonn\Ž $E$ et de l'ensemble totalement ordonn\Ž $\cal C(E)$ des coupures dans $E$. Soient $x\in E$ et $c = (A|B) \in \cal C(E)$; il n'y a que deux possiblit\Žs pour $x$ : ou bien $x$ appartient \ˆ la partie $A$, ou bien $x$ appartient \ˆ la partie $B$ car $(A | B)$ est une coupure. Si $x\in A$, on \Žcrit $x\leqs c$ et si $x\in B$, on \Žcrit $c \leqs x$.

\

\noi On a ainsi trois relation binaires d\Žsign\Žes par $\leqs$ : l'une sur $E$, une autre sur $\cal C(E)$ et une troisi\me entre les \Žl\Žments de l'un et l'autre ensemble. En les combinant toutes les trois,  on obtient une relation d'ordre totale sur l'ensemble $\cal A=E \cup \cal C(E)$ comme on peut le v\Žrifier simplement, pas \ˆ pas. [Il faut le faire soi-m\me !] Pour ainsi dire, on a un {\it amalgame} des deux ensembles totalement ordonn\Žs, $E$ et $\cal C(E)$. Pour l'ordre strict associ\Ž, on a $(A|B) < (C|D)$ si et seulement si l'on a  $A\inc C \et C\stm A \neq \vide$.

\

\lopar On verra, plus loin, comment la construction de la classe des  nombres de {\sc Conway} est une simple r\Žp\Žtition de ce proc\Žd\Ž d'amalgamation! Tu remarqueras que j'utilise un vocabulaire qui n'est pas toujours celui des usages courants. Ne t'inqui\te pas, cela restera entre nous\}

\

\su{Abr\Žviations} Soient  $X$ et $Y$ deux parties quelconques d'un ensemble totalement ordonn\Ž  et soit $z$ un des \Žl\Žments de cet ensemble. Lorsque l'on a $x < z$ pour chaque $x\in X$, on \Žcrit $X < z$, pour abr\Žger. De m\me, $z < Y$ est une abr\Žviation de ($z < y$ pour chaque $y\in Y$). On \Žcrit $X<Y$ pour dire que, pour tout $x\in X$ et tout $y\in Y$, on a $x<y$. De m\me, on \Žcrira $X\leqs Y$ lorsque l'on a  $x\leqs y$, pour tout $x\in X$ et tout $y\in Y$. 

\

\noi En particulier, lorsque $z = ( X|Y)$ est une coupure dans $E$, il est clair que l'a on $X < z < Y$, dans l'amalgame $\cal A=E \cup \cal C(E)$.

\

\su{Corps ordonn\Žs} Par d\Žfinition, {\bf un corps ordonn\Ž} est un corps commutatif $K$ muni d'une {\bf relation d'ordre total}, $\leqs$, compatible avec la structure de corps. Ce qui veut dire que, pour tous \Žl\Žments $x,y,z,$  de $K$, on a  :
$$x\leqs y \implies x+z \leqs y+z$$
$$x\geqs 0  \et y\geqs 0  \implies  xy\geqs 0.$$
Parmi les corps ordonn\Žs, on en distingue certains que {\sc Bourbaki} nomme {\bf les corps ordonn\Žs maximaux} et que l'on appelle aussi, couramment, {\bf les corps r\Žellement clos}. L'exemple type de corps ordonn\Ž  maximal est le corps $\R$ des nombres r\Žels. La th\Žorie bien connue des corps ordonn\Žs maximaux dit, en particulier, que tout corps ordonn\Ž se plonge (se prolonge) en un corps ordonn\Ž maximal. C'est ainsi  que le corps ordonn\Ž  des nombres rationnels $\Q$  est plong\Ž dans le sous-corps r\Žellement clos de $\R$ form\Ž des nombres alg\briques r\Žels.

\

\head{Notations, conventions, extensions}

\

\su{Intervalles} On distingue, en particulier,  deux types d'intervalles : l'intervalle ferm\Ž $[x,y]$ et l'intervalle ouvert $]x,y[$.  Par d\Žfinition, l'intervalle ouvert $]x,y[$ est l'ensemble des \Žl\Žments $z$ tels que $x < z < y$.  L'intervalle ferm\Ž $[x,y]$ est l'ensemble des \Žl\Žments $z$ tels que $x\leqs z \leqs y$. Le premier est contenu dans le second. Ils peuvent \tre vides ! {\it et}, par exemple, {\it l'intervalle $[x,y]$ l'est si et seulement si l'on a $y<x$}.

\

\noi{\bf Ces d\Žfinitions et notations s'\Žtendent aux classes totalement ordonn\Žes.} Pour deux ensembles de nombres $X$ et $Y$, on d\Žsigne par $]X,Y[$, resp., $[X,Y]$, l'ensemble des nombres $z$ tels que $X < z < Y$, resp. $X \leqs z \leqs Y$. De m\me pour deux classes de nombres, $\sX$ et $\sY$.

\

\noi {\bf Que peut-il se passer entre deux nombres, deux ensembles de nombres ou deux classes de nombres ?} On r\Žpondra \ˆ ces questions dans la suite du texte !

\

\noi J'utilise les notations classiques pour les ensembles classiques. En particulier :

$\N$ \ ensemble des entiers naturels \ , \ $\Z$ anneau des entiers relatifs \ , \ $\Q$ corps des nombres rationnels \ , \ $\R$ corps des nombres r\Žels.

\

\noi J'introduis les notations particuli\res suivantes :

$\sS$  \  classe des nombres \ , \  $\sM$ \ classe des noms \ , \ $\sO$ \ classe des ordinaux.

\

\lopar Ne t'\Žtonne pas : je d\Žsigne la classe des nombres  par $\sS$ car $\sN$ pourrait pr\ter \ˆ confusion et que, comme tu l'as  lu sans doute, {\sc Knuth} a appel\Ž {\bf surr\Žels} les nombres de {\sc Conway}.\}

\

\noi On a \Žtendu aux classes totalement ordonn\Žes les notations telles $\sX >x$ ou telles $\sX <\sY$, de mani\re naturelle. On \Žtend de m\me les notions de cofinalit\Ž et de co\•nitialit\Ž aux  classes totalement ordonn\Žes, comme suit.

\su{Parties cofinales et parties co\•nitiales} Soient $\sX$, $\sY$, deux sous-classes d'une m\me classe totalement ordonn\Že. On dit que $\sY$ est cofinale \ˆ $\sX$ lorsque, pour chaque $x\in \sX$, il existe un $y\in \sY$ tel que $x\leqs y$. On dit que les deux sous-classes, $\sX$ et $\sY$,  sont {\bf cofinales} lorsque, chacune est cofinale \ˆ l'autre [pour chaque $x\in \sX$, il existe un $y\in \sY$ tel que $x\leqs y$ et, vice versa, pour chaque $y\in \sY$,  un $x\in \sX$ tel que $y\leqs x$]. On dit que $\sX$ et  $\sY$ sont co\•nitiales lorsqu'elles sont cofinales pour  l'ordre total $\geqs$ inverse. Bien entendu, on dit que $\sY$ est co\•nitiale \ˆ $\sX$ lorsqu'elle lui est cofinale pour l'ordre inverse !

\noi On verra ci-dessous comment les classes $\sO$  et $\sS$ sont cofinales !

\

On construit la classe $\sS$ par r\Žcurrence, une r\Žcurrence transfinie.

\

\head{2 \ En route pour la r\Žcurrence}

\

Ici, $\alpha$ d\Žsigne un  ordinal,  $\alpha \in \sO$, bien entendu.

\

\noi  \`A l'\Žtape $\alpha$ de la r\Žcurrence,    on prend pour $S_\alpha$  l'ensemble des coupures dans l'ensemble  totalement ordonn\Ž $T_\alpha = \cup _{\beta<\alpha} S_\beta$, r\Žunion des ensembles $S_\beta$ construits aux  \Žtapes pr\Žc\Ždentes, $\beta < \alpha$. On obtient l'ensemble totalement ordonn\Ž $T_{\alpha + 1} = T_\alpha \cup S_\alpha$ comme amalgame, par le proc\Žd\Ž d\Žj\ˆ d\Žcrit dans le {\sc Prologue}.

\

\su{Explications et vocabulaire}  Dans le {\sc Prologue}, on a dit comment on ordonne $S_\alpha$, l'ensemble des coupures, de mani\re naturelle,  et comment on d\Žfinit  l'ordre total sur $T_\alpha \cup S_\alpha$ par amalgame !    Les nombres nouvellement cr\Ž\Žs  \ˆ l'\Žtape $\alpha$ constituent l'ensemble $S_\alpha$ : c'est la  {\bf g\Žn\Žration $\alpha$}. L'ensemble $T_\alpha$ est la r\Žunion des  g\Žn\Žrations pr\Žc\Ždentes, $S_\beta$ pour $\beta < \alpha$. 

\

\noi Pour un nombre $a$ de la g\Žn\Žration $\alpha$, on \Žcrira $g(a) = \alpha$, l'indicatif de sa g\Žn\Žration. De mani\re imag\Že, pour deux nombres , $a$ et $b$, on dira que $a$ est plus vieux que $b$ lorsque l'on  a $g(a)<g(b)$, qu'ils ont le m\me \‰ge si $g(a) = g(b)$.

\

\noi Autrement dit, les nombres qui ont un m\me \‰ge appartiennent \ˆ une m\me g\Žn\Žration $S_\alpha$, et les plus vieux constituent l'ensemble $T_\alpha$.

\

\noi Pour amorcer la r\Žcurrence, on part de l'ensemble vide, $\vide$, dont l'unique coupure est $(\vide|\vide)$ : c'est le premier des nombres, il est de g\Žn\Žration 0, autrement dit, $S_0 = \{(\vide|\vide)\}$. Il y a deux coupures dans l'ensemble $S_0$ : \ˆ savoir $(\vide|(\vide|\vide))$ et $((\vide|\vide)| \vide)$.  Ce sont les deux seuls nombres de g\Žn\Žration 1. Les trois premiers nombres  sont ainsi plac\Žs, d'apr\s les d\Žfinitions, dans l'ordre total suivant :
$$(\vide|(\vide|\vide)) < (\vide|\vide) < ((\vide|\vide)| \vide).$$
On voit bien que cela finira vite par devenir illisible, si l'on n'y prend garde ! Aussi, introduit-on des {\it abus d'\Žcriture} : au lieu de $(\vide|\vide)$, on \Žcrit simplement $( \ | \ )$ et on pose $0 = ( \ | \ )$,  autrement dit, on d\Žsigne le premier nombre par l'entier $0$. De m\me, au lieu de $(\vide|(\vide|\vide))$ et  de $((\vide|\vide)| \vide)$, on \Žcrit $( \ | \ 0)$ et $(0\ |\ )$  et l'on pose $-1 = ( \ | \ 0 ) \ , \  1 = (0 \ | \ )$. Les trois  nombres des g\Žn\Žrations 0 et 1  sont donc 
$$-1 < 0 < 1.$$
Ainsi de suite ... Allons jusqu'\ˆ la g\Žn\Žration $2$, sans commentaires.

\

\noi G\Žn\Žration 2 : $$-2 = ( \ |-1, 0, 1) < (-1 | 0, 1) <  (-1,0 |  1) < 2= (-1,0,1| \ ).$$
Voici l'ensemble des trois premi\res g\Žn\Žrations, 0, 1, 2, et son ordre total strict (respectant toutes les d\Žfinitions) :
$$-2 < -1 < ( \ |-1, 0, 1) < 0 < (-1 | 0, 1) < 1 < 2.$$
$$-2 < -1 < -1/2 < 0 < 1/2 < 1 < 2.$$
C'est \ˆ dessein que l'on d\Žsigne ces sept premiers nombres ainsi, en les identifiant aux nombres entiers  ou rationnels correspondants.   Comme on le verra par la suite,  plus bas,  les ensembles classiques $\Z$ et $\Q$ s'identifient \ˆ  leurs copies dans la classe $\sS$ des nombres. On verra en particulier, apr\s avoir d\Žfini l'addition, que l'on a bien $1+1 = 2$ et $1/2 + 1/2 = 1$. 

\

\head{3 \ Nombres oppos\Žs}  

\
 
\noi  On reprend les notations du paragraphe 2, en particulier, $T_\alpha = \cup _{\beta<\alpha} S_\beta$.

\

\noi L'oppos\Ž d'un nombre  $x$ est un nombre d\Žsign\Ž par $-x$. On d\Žfinit cette notion  par r\Žcurrence. On suppose $-x$ d\Žfini pour tout $x\in T_\alpha$. Pour chaque partie $X\inc T_\alpha$, on pose $-X = \{-x : x \in X\}$. Enfin, lorsque $x = (A | B)\in S_\alpha$ est un nombre,   on prend  $-x = (-B | -A)$, sachant que $(-B | -A)$ est une coupure dans $T_\alpha$. Bien entendu, on a $-0 = 0$ puisque $0 = ( \ | \ )$. Ainsi, l'oppos\Ž de $1 = (0 | \ )$ est $-1 = ( \ | 0)$; de m\me que l'oppos\Ž de $2 = (-1, 0, 1| \ )$ est  $-2 = ( \ |-1, 0,1)$, conform\Žment aux notations adopt\Žes ci-dessus. On voit assez clairement que,   pour tout nombre  $x$, on a $-(-x) = x$ : tout nombre est l'oppos\Ž de son oppos\Ž !

\

\su{Remarque} Chaque g\Žn\Žration $\alpha$ est un \emph{ensemble} de nombres. Parmi eux, il y en a un plus grand que tous les autres, $(T_\alpha | \ )$, de sorte que le nombre $-(T_\alpha | \ ) = ( \ | T_\alpha)$ est le plus petit de sa g\Žn\Žration.  On identifie le nombre  $(T_\alpha | \ )$ au nombre ordinal $\alpha$, de sorte que $\sO$, la classe  des  ordinaux, se pr\Žsente comme une  sous-classe de la classe $\sS$  des nombres. Comme on l'a mentionn\Ž plus haut, ces deux classes sont cofinales !

\

\noi L'ensemble $T_{\alpha + 1}$ des nombres de toutes les g\Žn\Žrations $\beta \leqs \alpha$ est ainsi une partie du segment $[-\alpha, \alpha]$.

\

\lopar Ne pas confondre {\it plus vieux} qui est relatif aux g\Žn\Žrations, avec les expressions {\it plus petit} ou {\it plus grand} qui sont relatives \ˆ l'ordre !\}

\

\head{4 \ Noms et synonymes}

\

Pour la commodi\Ž, j'introduis les notions de {\bf nom}  et de {\bf synonyme}. On s'en servira, en particulier pour la d\Žfinition de l'addition, de la multiplication et autres op\Žrations dans la classe $\sS$ des nombres !

\

Voici d'abord un r\Žsultat simple et utile.

\su{Lemme} {\sl Entre deux nombres $x < y$ d'une m\me g\Žn\Žration  $\alpha$ donn\Že, il y a toujours un nombre $z$ plus vieux, c'est \ˆ dire d'une g\Žn\Žration ant\Žrieure, $\beta < \alpha$.}

\su{D\Žmonstration} Soient $x = (A|B) < y = (C,D)$ deux nombres d'une m\me g\Žn\Žration $\alpha$. On a $A\inc C \inc T_\alpha$ et $A \neq C$, par d\Žfinition. Il existe donc au moins un nombre $z \in C\stm A\inc T_\alpha$ :  ce  $z$  est donc d'une g\Žn\Žration $\beta$ pr\Žc\Ždant $\alpha$, autrement dit, $\beta < \alpha$, et $z\in B$ de sorte que l'on a $x < z < y$, d'apr\s les d\Žfinitions !\qed

\

\noi On \Žtend ce r\Žsultat aux cas des ensembles  de nombres, $X,Y$, parties  de la classe  totalement ordonn\Že $\sS$. On a d\Žfini {\bf l'intervalle} $]X,Y[$ ayant pour extr\Žmit\Žs $X,Y$, comme \Žtant la classe des nombres $z$ tels que l'on ait $X <  z < Y$. Bien entendu, cet intervalle peut, parfois, \tre vide.

\

\su{Th\Žor\me du nombre m\Ždiateur} {\sl Soient $X,Y$, deux ensembles de nombres tels que $X < Y$. Il existe au moins un nombre  $z$ pour lequel on a $X<z<Y$ et, parmi tous ces nombres, il en existe un, et un seul,  plus vieux que les autres.} On dira que c'est {\bf le nombre m\Ždiateur} entre $X$ et $Y$.

\

\su{D\Žmonstration \emph{d\Žtaill\Že}} On se souvient que  la classe $\sO$ des nombres ordinaux est \emph{bien ordonn\Že} et que $g(z)$ d\Žsigne le num\Žro de la g\Žn\Žration du nombre $z$. Les $g(z)$, pour $z\in X\cup Y$,  forment un  \emph{ensemble} $G$ de nombres ordinaux puisque $X\cup Y$ est un \emph{ensemble}, de sorte qu'il existe des ordinaux $\gamma$ tels que $\gamma > G$ et  parmi ceux-l\ˆ un plus petit que tous les autres ! 

Soit $\alpha$ le plus petit ordinal tel que l'on ait $g(z)<\alpha$ pour tout $z\in X\cup Y$. 

On a ainsi $X \inc T_\alpha \et Y\inc T_\alpha$ 
et il existe un couple (parfois m\me plusieurs) de parties $A,B$,  de $T_\alpha$ tels que  l'on ait 
$$X \inc A \ , \ Y\inc B \ , \ A\cap B = \vide \ , \ A\cup B = T_\alpha,$$ 
de sorte que $z = (A|B)$ est une coupure dans l'ensemble $T_\alpha$, autrement dit $z$ est un nombre de g\Žn\Žration $\alpha$ et l'on a ainsi $X <  z < Y$, d'apr\s les d\Žfinitions !  Parmi ces $z$, il y en a un seul  plus vieux que tous les autres, d'apr\s le lemme ci-dessus car, s'il y en avait deux distincts ils seraient de la m\me g\Žn\Žration ! \qed

\

\lopar J'ai pr\Žf\Žr\Ž l'appellation {\it nombre m\Ždiateur} \ˆ celle de {\it nombre interm\Ždiaire} ! On pourrait lui donner tout autre nom qui plairait davantage !\}

\su{Les noms}  On dit qu'un   couple donn\Ž, $X, Y$,  quelconque, d'ensembles de nombres  est un {\bf nom} lorsque que l'on a $X < Y$. Pour le noter, on  utilisera le symbole $\{X | Y\}$. La notion et la notation sont semblables \ˆ celles des coupures, mais diff\Žrentes ! Insistons bien !

\

\noi Pour chaque nom $\{X | Y\}$, il existe toujours au moins un nombre $z$ entre $X$ et $Y$, autrement dit tel que $X < z < Y$ et, parmi tous ces nombres, il en existe un, et un seul, plus vieux que tous les autres, le nombre m\Ždiateur :  c'est le th\Žor\me du nombre m\Ždiateur.

\

\noi {\bf On dira que  $\{X | Y\}$ \emph{est un des noms} de ce nombre m\Ždiateur.}

\

\noi D'autre part, pour chaque nombre  $a = (A|B)$, le nom $\{A|B\}$  est  un nom de $a$  : on dira que c'est le nom {\bf intime} du nombre $a$ lui-m\me. 

\

\noi Un nombre  peut avoir aussi plusieurs autres  noms  !

\

\su{R\Žsumons} Soit $\{X|Y\}$ un nom : on a $X < Y$ par d\Žfinition.  
Le nombre  $z = (A|B)$ a pour nom $\{X | Y\}$ veut dire que $z$ est le nombre m\Ždiateur entre $X$ et $Y$.  Pour cela, il faut et il suffit que les trois conditions suivantes soit satisfaites :

\

\noi   on a $X < z < Y$,

  \noi pour chaque $t\in A$, il existe un $x\in X$ au moins tel que $x \geqs  t$, 
  
  \noi pour chaque  $t\in B$, il existe un $y\in Y$ au moins tel que $y \leqs t$.
  
  \

\su{On le v\Žrifie simplement} [Il faut le refaire soi-m\me !] Le nombre m\Ždiateur $z$ est le plus vieux parmi tous les nombres $u$ pour lesquels on a  $X < u <Y$. Autrement dit, pour chaque $t\in A\cup B$, on a ou bien non($X < t$) ou bien non($t <Y$). Pour $t\in A$, on a d\Žj\ˆ
$t < z< Y$, il faut donc, et il suffit, que l'on ait non($X < t$). De m\me, pour $t\in B$, il faut, et il suffit que l'on ait non($t < Y$). Il faut donc et il suffit que, pour $t\in A$ on ait non($X < t$) et, pour $t\in B$ on ait non($t< Y$). Or, non($X< t$) \Žquivaut \ˆ $(\exists x\in X)(x\geqs t)$, et non($t < Y$) \Žquivaut \ˆ $(\exists y\in Y)(y \leqs t)$. \hfill{\bf Fin de la v\Žrification}

\

\su{Une condition n\Žcessaire et suffisante} En utilisant les notions de cofinalit\Ž et de co\•ntialit\Ž, on peut \Žnoncer cette caract\Žrisation sous la forme suivante. Le nom $\{X | Y\}$ d\Žsigne le nombre $a = (A | B )$ si et seulement si :
$$\text{on a } X < a < Y \ , \ X \ \text{est cofinal \ˆ} \ A \ , \  Y \ \text {est co\•nitial \ˆ} \ B.$$
\su{Un crit\re suffisant} La condition suivante est suffisante pour qu'un nom $\{X | Y\}$ soit celui du nombre  $a = (A | B )$ :
$$X \et A \ \text{sont cofinaux} \ , \  Y \et B \ \text{sont co\•nitiaux}.$$
En effet, la condition ($A$ est cofinal \ˆ $X$) entra\"ne ($X < a$) et la condition ($B$ est co\•nitial \ˆ $Y$) entra\"ne ($a < Y$).

\

\su{Les synonymes} On dira que deux noms, $\{A|B\}$ et $\{C|D\}$,  sont synonymes lorsqu'ils d\Žsignent, tous deux,  le m\me nombre $a$ et on \Žcrira 
$$\{A|B\}\sim \{C|D\}\sim a.$$
\noi On identifie chaque nombre \ˆ son nom intime, de sorte que la classe des nombres, $\sS$,  se pr\Žsente comme une {\bf sous-classe} de la classe $\sM$ de tous les noms.   La relation $\sim$ est une relation d'\Žquivalence d\Žfinie sur la classe $\sM$. Cette relation d'\Žquivalence induit l'identit\Ž sur la sous-classe $\sS$  des nombres. Dans la classe des nombres, la synonymie n'est autre que l'identit\Ž ! Deux nombres  sont  synonymes si et seulement s'ils sont \Žgaux !!!

\

\su{Une condition suffisante pour la synonymie} Lorsque les deux parties, $ X,U$, sont cofinales et les deux parties, $Y,V$, sont co\•nitiales, les deux noms $\{X | Y\}$ et $\{U | V \}$ sont synonymes.

\

\noi Cela d\Žcoule du crit\re suffisant signal\Ž plus haut.

\

\su{Exemple} Soit $\{X|Y\}$ un nom. On se donne $u\in X$ et $v\in Y$ et on pose 
$$U = \{x \in X : x \geqs u \} \ , \ V = \{y\in Y : y \leqs v\}.$$
Alors $\{U|V\}$ est synonyme de $\{X|Y\}$.

\

\su{Exercice} Soit $\{X|Y\}$ un nom du nombre $a = (A|B)$. Retrouver les parties $A$ et $B$ \ˆ partir des parties $X$ et $Y$. Autrement dit, reconstituer le nom intime d'un nombre \ˆ partir de l'un quelconque de ses noms.

\

\head{5 \ L'addition}

\

\lopar On le sait, il est d'usage d'\Žtendre une op\Žration binaire $\bot$ donn\Že  sur un ensemble $E$ quelconque, de mani\re naturelle, en une op\Žration sur les parties de cet ensemble, en posant :
$$X \ \bot \ a = \{x \ \bot \ a : x\in X\} \ , \ X \ \bot \ Y = \{x \ \bot\ y : x\in X, y\in Y\}.$$
 C'est le cas, par exemple, pour les groupes, les anneaux ou les espaces vectoriels !\}

 \

\noi On d\Žfinit une addition, $a + b$, dans la classe $\sS$   des nombres  en associant, \ˆ chaque couple de nombres, $(a,b),$ un nombre  $c$ et l'on \Žcrit $a + b = c$.  La d\Žfinition se fait par r\Žcurrence. Une r\Žcurrence transfinie,  laborieuse, quelle que soit la mani\re dont on l'aborde. Elle n\Žcessite de tr\s nombreuses v\Žrifications, fastidieuses,  mais il vaut mieux les faire soi-m\me, une fois au moins dans sa vie !

\

\noi \`A l'\Žtape $\alpha$ de la r\Žcurrence, l'addition est d\Žj\ˆ d\Žfinie pour tous les nombres de $T_\alpha$. On poursuit alors en d\Žfinissant l'addition dans  l'ensemble $T_{\alpha + 1} = T_\alpha\cup S_\alpha$.

\

\noi Voici le c\oe ur de la construction.

\

\noi Soient  $\{A | X \}$ le nom intime de $a$ et $\{B | Y \}$ celui de $b$. On d\Žfinit la somme $a + b$ comme \Žtant le nombre dont  $\{A+B | X + Y\}$  est un nom. Autrement dit, $a + b = (A|X) + (B|Y) \sim \{A+B | X + Y\}$.

\

\noi Le d\Žpart de la r\Žcurrence se fait avec la premi\re g\Žn\Žration, $\{0\}$.  On part de  $0 = (\ | \ )$ de sorte  que $0 + 0$ a pour nom  $\{ \  | \ \}$ lequel est un nom de $0$, d'o\ $0 + 0 = 0$ !

\

\su{Explication} On veut que l'addition satisfasse la condition suivante :
$$(a \leqs x \et b \leqs y) \implies (a+b \leqs x + y)\ !$$
Cela \Žclaire la d\Žfinition que l'on en donne : cette d\Žfinition est, pour ainsi dire, {\it n\Žcessaire et suffisante} !

\

\noi Cette d\Žfinition permet de voir, sans grand frais, que  $0+ a= a = a+0$, pour tout nombre  $a$. On voit \Žgalement que $a+ (-a)$ a pour nom $\{A-B | B-A \}$ qui est un nom de $0$ !

\noi  Plus g\Žn\Žralement, quels que soient les noms choisis $\{X|U\}$  et $\{Y|V\}$ pour $a$ et $b$ respectivement,  aura   : 
$$a+b \sim \{X+b, Y+a | U+b, V+a\}.$$
On montre que cette addition est commutative, associative, qu'elle est compatible avec l'ordre total des nombres, $0$  est un \Žl\Žment neutre et l'on a $x + (-x) = 0$.  Autrement dit, munie de cette op\Žration,  la classe des nombres poss\de la  structure d'un {\bf groupe totalement ordonn\Ž},
commutatif.

\

\noi Petit exemple : $\{-1,0 | 1\} + \{-1,0 | 1\} = \{-2,-1, 0 | 2\}$ qui est un nom du nombre $(0 | \ ) = 1$. Cela justifie les notations utilis\Žes  ci-dessus, au paragraphe 2, pour d\Žsigner les nombres des trois premi\res g\Žn\Žrations !

\newpage

\

\head{6 \ La multiplication}

\

Par un proc\Žd\Ž semblable au pr\Žc\Ždent, on d\Žfinit une multiplication entre les nombres en associant, \ˆ chaque couple de nombres , $(a,b),$ un nombre  $c$ et l'on \Žcrit $a.b = c$. On montre que cette op\Žration est commutative, associative, distributive par rapport \ˆ l'addition [autrement dit, $x.(y+z) = x.y + x.z$], que tout nombre  $x\neq 0$ poss\de un inverse $1/x$ [autrement dit, $x.1/x = 1$]  et poss\de la propri\Žt\Ž suivante :
$$0 >0  \et  y > 0 \implies x.y > 0, \pourtous x,y.$$
{\bf En r\Žsum\Ž, cela revient \ˆ dire que la classe des nombres munie de l'ordre total, de l'addition et de la mutiplication ainsi d\Žfinies v\Žrifie les axiomes des corps ordonn\Žs } :  la d\Žfinition est donn\Že, ci-dessus, dans le {\sc Prologue}.

\

\su{Pr\Žparation} \`A certains  nombres dans la classe $\sS$, on peut donner des noms particuliers, carat\Žristiques. Par exemple, soit $a= (A|B)$ un nombre strictement positif, $a> 0$, d'une g\Žn\Žration $\alpha$ donn\Že. En prenant $X = A\cap \ ]0, a[$, on a $B > a > X >0$ o\ $X$ et $A$ sont cofinaux [comme on peut le v\Žrifier], de sorte que $\{X | B \}$ est encore un nom de $a$. 
 Tout nombre strictement positif $a > 0$ poss\de ainsi des noms $\{X | Y\}$ pour lesquels on a $Y> a > X>0$. Ce sont des noms {\bf caract\Žristiques} des nombres strictement positifs. \emph{Cum grano salis}, on pourra les appeler des {\bf dextronomes}.

\

\su{Esquisse du proc\Žd\Ž}Cette fois, j'omets beaucoup de d\Žtails et les v\Žrifications longues et fastidieuses, pourtant sans grand d\Žtour. Il vaut mieux que l'on entre soi-m\me dans ces d\Žtails pour les d\Žm\ler et faire les v\Žrifications n\Žcessaires.   Voici quelques indications.

\

\noi On  commence par d\Žfinir $0.x = x.0 = 0$, pour tout nombre $x$. Puis on d\Žfinit le produit  pour les nombres strictement positifs, par r\Žcurrence.  On compl\te enfin la d\Žfinition en  prenant $x.(-y) = -x.y$.  On suppose le produit $x.y$ d\Žj\ˆ d\Žfini pour tous $x >0$  et $y > 0$  pris dans $T_\alpha$, avec les propri\Žt\Žs voulues. Pour toutes parties $X$ et $Y$ de $T_\alpha$ telles que $X >0$ et $Y > 0$, on pose :
$$X.Y = \{x.y : x\in X, y\in Y\}.$$
C'est un ensemble de nombres strictement positifs. 
On utilise des dextronomes $\{X | U\ \}$ de $a$ et $\{Y | V \}$ de $b$. 
On d\Žfinit le produit $a.b$ comme \Žtant le nombre  qui a pour nom $\{X.Y| U.V \}$.  Autrement dit,
$$a.b \sim \{X.Y| U.V \}.$$
On  v\Žrifie que ce produit d\Žpend seulement de $a$ et de $b$, pas du choix de leurs dextronomes. 

\

\noi Comme on peut le v\Žrifier \Žgalement, l'inverse $1/a$ de $a > 0$, poss\de un nom $\{E | F \}$ o\ $E = \{z >0 : z.A < 1\}$ et $F = \{z > 0 : z.B > 1\}$.  

\

 \noi Ainsi de suite ... J'omets le reste. \hfill{\bf Fin de l'esquisse.}

\

\head{7 \ G\Žn\Žalogie}

\

Tu trouve, ci-dessous, une copie de la page 11 du livre de {\sc Conway}, {\sl On numbers and games} (la r\Žf\Žrence se trouve plus bas, \ˆ la fin de ma lettre). C'est un dessin qui repr\Žsente le haut de {\bf l'arbre g\Žn\Žalogique} des nombres selon {\sc Conway}.  

\includegraphics[width=5.3in]{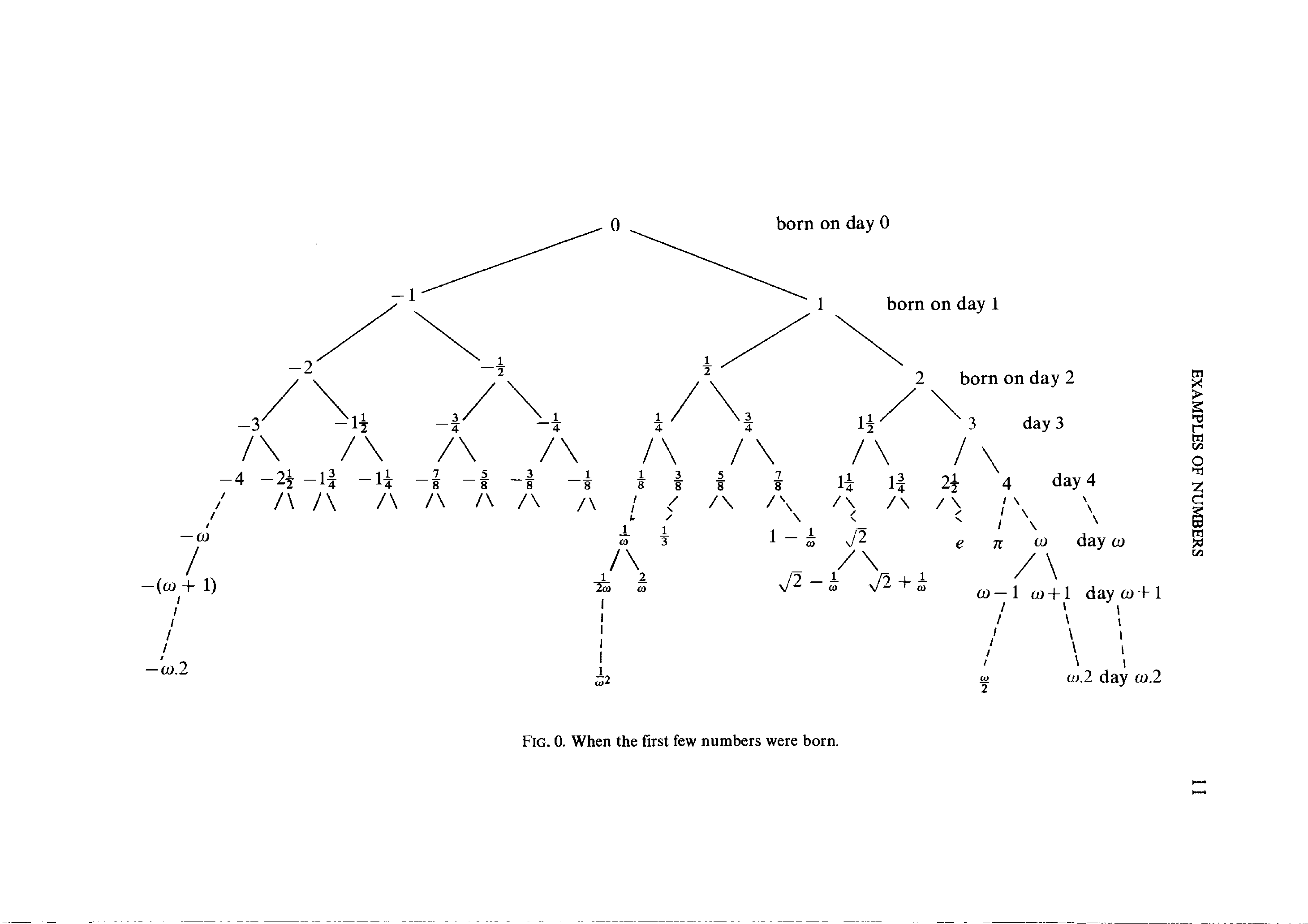}

\

\noi On y voit, en particulier, comment un nombre donn\Ž quelconque,  d'une g\Žn\Žration $\alpha$  poss\de {\bf deux successeurs} dans la g\Žn\Žration   $\alpha +1$ suivante    : l'un \ˆ droite et l'autre \ˆ gauche. Par exemple, les deux successeurs du nombre $0$ sont $+1$ et $-1$. De m\me, ceux de $3$ sont $4$ et $\frac{5}{2}$, lequel est ici d\Žsign\Ž $2\frac{1}{2}$ \ˆ l'anglo-saxonne !

\

\lopar Pour Conway chaque nombre poss\de un \emph{birthday}, le jour o\ \emph{il est n\Ž}. Conway utilise ainsi l'expression $x$ {\it was born on day} $\alpha$ l\ˆ ou je dis $x$ est de la g\Žn\Žration $\alpha$.\}

\

\noi Voici quelques explications suppl\Žmentaires.

\

\noi Soit $x= (A|B)$  un nombre d'une g\Žn\Žration donn\Že $\gamma > 0$. Pour chaque ordinal $\alpha < \gamma$, posons $A_\alpha = T_\alpha \cap A$ et $B_\alpha = T_\alpha \cap B$. Ainsi $x_\alpha = (A_\alpha | B_\alpha)$  est une coupure dans $T_\alpha$ : c'est un nombre  de la g\Žn\Žration $\alpha$.  Cela est clair. Pour $\alpha = 0$, on a $x_0 = 0$, bien entendu.

\

\noi On obtient la suite $s_\gamma(x)= (x_\alpha)_{0\leqs \alpha < \gamma}$ de nombres, index\Že par les ordinaux $ \alpha < \gamma$. Il est commode de dire que, pour un $\beta$ donn\Ž, $x_\beta$ est {\bf l'ascendant} de $x$ dans la g\Žn\Žration $\beta$,  et que $x$ est {\bf l'un des descendants} de $x_\beta$ dans la g\Žn\Žration $\alpha$.

\

Il convient  d'appeler la suite  $s_\gamma(x)$ {\bf l'ascendance} de $x$ ou encore {\bf la liste des ascendants} de $x$. Il est clair que l'ascendance de $x_\alpha$ n'est autre que la suite $(x_\beta)_{\beta < \alpha}$. Le descendant $x_{\alpha + 1}$ de $x_\alpha$ est, pour ainsi dire, {\it engendr\Ž} par $x_{\alpha}$.

\su{La relation p\re$|$enfant} Soient $x=(A|B)$ et $y = (C|D)$ deux nombres  des g\Žn\Žrations $\alpha$ et $\alpha + 1$, respectivement. Si $x$ est un ascendant de $y$, on peut dire, de mani\re imag\Že,  que $x$ est {\bf le p\re} de $y$, ou sa m\re, si l'on y tient ! De m\me, lorsque $y$ est un descendant de $x$,  on dira que $y$ est {\bf l'un des enfants} de $x$. Pr\Žcisons davantage la relation {\bf p\re$|$enfant}.

\

\noi Ne l'oublions pas, $x = (A|B)$ est une coupure dans  $T_\alpha$ et   $x$ appartient ainsi \ˆ $S_\alpha$. D'autre part, $y = (C | D)$  est une coupure dans l'ensemble $T_{\alpha + 1} = T_\alpha\cup S_\alpha$. Or, $y$ est un enfant de $x$ si et seulement si la coupure $(A|B)$ est la trace  de la coupure $(C|D)$ sur $T_\alpha$, par d\Žfiniton ! S'il en est ainsi, puisque $x$ appartient \ˆ $S_\alpha\inc T_{\alpha + 1}$, il n'y a que deux possibilit\Ž : on bien $x$ est le plus grand \Žl\Žment de l'ensemble $C$ et on aura $x < y$;  ou bien il est le plus petit \Žl\Žment de l'ensemble $D$ et on aura $y < x$. 

\

\noi Ainsi, un nombre $x$ quelconque engendre toujours deux enfants, et seulement deux,  disons $x_-$ et $x_+$, et l'on a  $x_- < x < x_+$.  Pour $x = 0$, ces deux enfants sont $-1$ et $+1$.

\

\noi Dans cet arbre, on appelle {\bf lign\Že} toute une suite de nombres  $x_\alpha$ index\Žs par les ordinaux $\alpha \geqs  0$,  o\ $x_{\alpha + 1}$ est l'un des deux descendants directs de $x_\alpha$.

\

\noi Par exemple, sur le c\™t\Ž droit de l'arbre, figure la lign\Že des nombres ordinaux,  $0, 1, \dots, \omega, \omega + 1, \dots, \omega. 2, \dots$. 

\

\noi Soient $\sA$   et $\sB$ deux sous-classes de la classe $\sS$ \ des nombres telles que  l'on ait  $\sA < \sB$  \ et $\sS = \sA\cup \sB$.  C'est une coupure dans la classe $\sS$,  une coupure \Žnorme, une sorte de faille ! Cette faille d\Žlimite une sorte de lign\Že transfinie $(x_\alpha)_{\alpha \in \sO}$   o\  $x_\alpha = (A_\alpha | B_\alpha) \ , \  A_\alpha = T_\alpha \cap \sA$ et $B_\alpha = T_\alpha \cap\sB$. 

\

\su{Exercice} Chercher la faille qui d\Žlimite la lign\Že des ordinaux !!! 

\

\head{8 \ Les g\Žn\Žrations finies}

\

Ce sont les premi\res g\Žn\Žrations, $S_0, S_1, S_2, \dots, S_n, \dots$, ind\Žx\Žes par les entiers naturels, autrement dit, les ordinaux finis $0,1,2,\dots, n,\dots$. La g\Žn\Žration $S_n$ est une partie de l'intervalle $[-n,n]$.  On a $|S_n| = 2^n$. Autrement dit, il y a $1,2,4,\dots,2^n, \dots$ nombres, respectivement, dans les g\Žn\Žrations $S_0, S_1, S_2, \dots, S_n, \cdots$. Leur r\Žunion est $T_\omega = \cup_{n\in \N} S_n$. On a $\Z\inc T_\omega$. Plus pr\Žcis\Žment, $T_\omega$ est l'anneau des nombres dyadiques :
$$T_\omega = \left\{\frac{k}{2^h} : k\in \Z , h\in \N\right\}= \D.$$
Cela se v\Žrifie simplement, sans d\Žtour, sachant que $\Z$ est une partie de $T_\omega$. Si $a\in T_\omega$, sa moiti\Ž $a/2$ appartient aussi \ˆ $T_\omega$.

\

\su{Les nombres r\Žels} Pour les besoins de la cause, on distingue parmi les nombres dyadiques, l'ensemble suivant :
$$\mathbf D = \{1/2^n : n \in \N \}.$$
On dit que $a$ est un nombre {\bf r\Žel} lorsqu'il existe un entier naturel $n$ tel que $-n \leqs a\leqs n$ et que le nom suivant est celui de $a$ :
$$\{a-\mathbf D | a+\mathbf D \} \sim a.$$
On montre que l'ensemble des nombres r\Žels ainsi d\Žfini est une copie dans $\sS$ du corps $\R$ des nombres r\Žels, classique !

\

\

\head{9 \ Forme normale}

\

Tout nombre peut \tre repr\Žsent\Ž sous la forme d'une s\Žrie transfinie, sa {\bf forme normale}. Cela rend ces nombres plus \emph{intelligibles}. Pour cela, on introduit les notions suivantes.

\

\su{Les ordres de grandeur} Comme dans tout corps ordonn\Ž,  on d\Žfinit la valeur absolue $|a|$ de $a$ comme suit :
$$|a| = a \si a\geqs 0 \ , \ |a| = -a \si a\leqs 0.$$
Pour les questions de comparaison des nombres, on s'en tient aux valeurs absolues. On peut ainsi se limiter \ˆ la classe $\sS_+ = \{a : a\geqs 0\}$ des nombres positifs.

\

\noi Soient $a, b$, deux nombres dans $\sS_+$. On dit que $a$ et $b$ sont {\bf commensurables} lorsqu'il existe un entier $n$ tel que $a \leqs nb$ et un entier $m$ tel que $b \leqs ma$. C'est une relation d'\Žquivalence. Les classes d'\Žquivalence correspondantes sont appel\Žes {\bf les ordres de grandeur} du corps ordonn\Ž. On dit que $a$ et $b$ sont de m\me ordre de grandeur lorsqu'ils sont commensurables. On \Žcrit $a << b$ lorsque, pour tout $n\in \N$, on a $na < b$, et on dit que $b$ est d'un ordre de grandeur sup\Žrieur \ˆ celui de $a$. Ainsi, entre deux nombres $a$ et $b$, une et une seule des relations suivantes tient : ou bien $a << b$ ou bien $b << a$ ou bien $a$ et $b$ sont commensurables.

\

\noi Les ordres de grandeur sont des intervalles, dans le sens suivant : si $a < b$ sont commensurables, tout l'intervalle $[a,b]$ est contenu dans le m\me ordre de grandeur que $a$ et $b$ !

\

\su{Exponentielle} On associe, \ˆ chaque nombre $a$, un nombre $\omega^a$, son {\bf exponentielle}. Cela se fait par r\Žcurrence. Comme d'habitude, pour chaque ensemble $X$ de nombres, on pose
$$\omega^X = \{\omega^x : x \in X\}.$$
Pour  un nombre $a = (X|Y)$ de g\Žn\Žration $\alpha$, on pose 
$$E = \{0\}\cup \N. \omega^X \ , \ F = \mathbf D.\omega^Y.$$
On prend pour $\omega^a$ le nombre qui a pour nom $\{E|F\}$ :
$$\omega^a\sim\{E| F\}.$$
On montre que l'on a

$\omega^a > 0$

$\omega^0 = 1$

$\omega^{a+b} = \omega^a . \omega^b$

$\omega^{-a} = 1/\omega^a$

$a < b \iff \omega^a << \omega^b$

Pour tout ordinal $\alpha \in \sO$, l'exponentielle $\omega^\alpha$ co\•ncide avec le nombre 

ordinal $\omega^\alpha$ .

\

\noi On \Žtablit aussi les deux r\Žsultats suivants. J'omets les d\Žmonstrations.

\

\su{Lemme}{\sl Pour chaque nombre $x > 0$, il existe un seul nombre $y$ tel que $x$ soit commensurable \ˆ $\omega^y$. De plus, $\omega^y$ est plus vieux ou aussi \‰g\Ž que $x$}.

\

\noi Plus pr\Žcis\Žment, le nombre $\omega^y$ est le plus vieux parmi tous les nombres qui ont le m\me ordre de grandeur que $x$.

\

\su{Th\Žor\me} {\sl Tout nombre $a$ poss\de une repr\Žsentation sous la forme d'une somme formelle 
$$a = \sum_{\alpha < \beta} \omega^{y_\alpha}. r_\alpha,$$
o\ $\beta$ est un ordinal, $(y_\alpha)_{\alpha < \beta}$ une suite strictement d\Žcroissante de nombres, et $(r_\alpha)_{\alpha < \beta}$ une suite de nombres r\Žels non nuls}. 

\

Cette repr\Žsentation, appel\Že {\bf forme normale}, est unique. Lorsque le nombre $a$ est un ordinal, cette repr\Žsentation co\•ncide avec la forme normale de Cantor pour les ordinaux !

\

\noi On pourrait fonder toute la th\Žorie des nombres de Conway sur cette repr\Žsentation. Ce ne serait pas la mani\re la moins attrayante.

\

\noi Les trois ordres de grandeur les plus familiers aux analystes sont
$$\omega^{-1} << \omega^{0} = 1 << \omega.$$
$$0 < \omega^{-1} < 1/n \leqs 1 \leqs n < \omega,\pourtout n\in \N^*.$$
En effet, $\omega ^{0} = 1$ est l'ordre de grandeur des nombres r\Žels. Tandis que le nombre $\omega^{-1}$ est un infiniment petit du premier ordre et $\omega$ est le premier entier infiniment grand. C'est un peu le langage que l'on utilise en analyse nonstandard de Robinson, c'est \ˆ dire la th\Žorie r\Žnov\Že des infinit\Žsmaux.

\

\noi L'un des points culminants de la th\Žorie des nombres de Conway est le r\Žsultat d\Žj\ˆ signal\Ž ci-dessus selon lequel la classe $\sS$ des nombres, [une classe propre, qui n'est pas un ensemble] poss\de les propri\Žt\Žs d'un corps ordonn\Ž r\Žellement clos lequel renferme une copie de chacun des corps ordonn\Žs qui sont des ensembles.

\

\noi Plut\™t que d'aller voir la d\Žmonstration dans l'un des nombreux textes o\ elle figure, on peut s'essayer \ˆ la faire soi-m\me ! R\Žussir serait la meilleure preuve que l'on a bien compris !!!

\

\head{Epilogue}

\

Pour en savoir davantage, on peut se reporter aux deux ouvrages suivants faciles \ˆ trouver.

\

\noi {\sc Norman L. ALLING}, {Foundations of analysis over surreal number fields}, Mathematics Studies 141, xvi + 373 pp., North Holland, 1987.

\

\noi {\sc J. H. CONWAY}, {On numbers and games}, ix + 230 pp.,
Academic Press Inc. 1976, reprinted 1979.

\

\noi La notion introduite par Conway, ses nombres, a re\c cu un tr\s bon accueil, un peu dithyrambique, parfois. Cependant, depuis longtemps on sait que \emph{tout ce qui est excessif est insignifiant}, comme le disait Talleyrand.  

\

\noi \`A l'\Žpoque,  Conway \Žtait d\Žj\ˆ bien connu, et appr\Žci\Ž. C'est Donald Knuth qui a donn\Ž le nom de nombres surr\Žels pour d\Žsigner les nombres de Conway. Il a \Žcrit une petite fable, une romance math\Žmatique, racontant l'histoire  de deux jeunes gens, Alice et Bill, qui d\Žcouvrent {\bf la pierre} grav\Že [une sorte de pierre de Rosette] qui m\ne \ˆ la cr\Žation de ces nombres. Cette histoire est bien connue. 

\

\noi Conway introduisait ses nombres avec un minimum d'axiomes et de d\Žfinitions que les contemporains avaient beaucoup appr\Žci\Žs. Pour l'illustrer, on trouvera dans les  {\sc Annexes} ci-dessous, trois pages de son livre qui semblent tout r\Žsumer. La plupart de ses ex\Žg\tes lui  ont emboit\Ž le pas et de nombreux auteurs, pour introduire les nombres de Conway, se servent encore de cette mani\re cursive et pas tr\s ais\Že, presque sibylline !

\

\noi Mais l'histoire ne s'arr\te pas l\ˆ. On apprend, dans le livre de Alling, que Conway a un pr\Žcurseur, Norberto Cuesta (1907-1989). En effet, dans un article de 1954, ce dernier, grand admirateur de Sierpinski, construisait d\Žj\ˆ les nombres de Conway, avant la lettre, dans l'article suivant :

\

\noi {\sc N. CUESTA}, {\sl Algebra ordinal}, Revista de la Real Academia De Ciencias Exactas, Fisicas Y Naturales, {\bf 58} \no 2 (1954) 103-145.

\

\noi C'est Charles Helou, ami et coll\gue qui a fini par  d\Žnicher ce texte pr\Žcieux, tr\s difficile \ˆ retrouver, avec l'aide de la  biblioth\que centrale de Penn State University. Je lui en suis infiniment reconnaissant. Ce texte est \Žcrit en espagnol, en Espagne, du temps de Franco ! J'en  ai fait la traduction en fran\c cais, avec l'aide de mon \Žpouse Claude.

\

\noi Nous voil\ˆ devant un m\me objet math\Žmatique, observ\Ž sous deux angles compl\Žmentaires, \ˆ une trentaine d'ann\Žes de distance ! On dirait que l'un porte son regard de bas en haut pour observer l'\Ždifice, et l'autre veut l'observer de haut en bas, en renversant la perspective ! 

\

\noi De mani\re cursive, je dirai que  l'approche de {\sc Conway} est synth\Žtique, celle de {\sc Cuesta} analytique !

\

\noi Qu'a bien pu dire {\sc Conway} en apprenant cette co\•ncidence ! Qu'a bien pu dire {\sc Cuesta Dutari} lui-m\me ! On aimerait bien le savoir. Mais on ne le  saura probablement jamais ! Il y a peut-\tre eu des propos \Žchang\Žs ou not\Žs, quelque part, on ne sait o\. Pourvu que la poussi\re des si\Žcles ne les recouvrent pas !

\

\noi Par un curieux hasard, ou une \Žtrange co\•ncidence, je viens de faire la connaissance, \ˆ travers la toile, d'un coll\gue math\Žmaticien, {\sc Ricardo P\Žrez-Marco},  qui a pu serrer la main de {\sc Cuesta} dans les rues de Salamanque, quand il \Žtait jeune. Je trouve cela \Žmouvant  !

\newpage

\noi Cher Ami, 

\

\noi Je sais d\Žj\ˆ que tu ne liras pas cette lettre de bout en bout; je devine que non !  Il n'emp\che, j'y ai mis l'essentiel de ce qu'il faut savoir pour comprendre ce que sont les nombres surr\Žels de {\sc Conway},  et c'est long ! 

\

\noi Je te rappelle le c\Žl\bre adage de {\sc Erd\šs} : \guil Tout le monde \Žcrit, personne ne lit !\guir

\

\noi Pourtant, si long soit-il, comme l'Iliade ou l'Odyss\Že, je crois que l'on peut lire un beau po\me jusqu'au bout, et m\me plusieurs fois de suite ! Mais qui a jamais vraiment d\Žj\ˆ \Žt\Ž jusqu'au bout d'un texte math\Žmatique ... C'est tr\s rare ! Il y faut une attention soutenue, une grande patience et une envie irr\Žpressible, irr\Žsistible ! Jamais, pour ainsi dire ...

\

\noi Je t'ai aussi pr\Žpar\Ž une liste de lectures, si jamais tu en as l'envie.

\

\head {Pour en savoir plus}

\

Voici  d'abord deux livres importants sur le sujet :

\

{\sc John Horton Conway}, {\sl On numbers and games}, Academic Press,
(1976), reprinted with corrections 1977, reprinted 1979, ix + 238 pp.

\

{\sc Norman L. Alling}, {\sl Foundations of analysis over surreal number fields},  North-Holland Mathematics Studies, (1987) {\bf 141}, xvi + 373 pp.

\

\noi Puis une liste plus fournie.

\

\noi[ 1 ] {\sc Norman L. ALLING}, {\sl Conway's field of surreal numbers},Trans.A.M.S., {287} \no 1 (1985) 365-386.

\

\noi [ 2 ] {\sc Norman L. ALLING}, {Foundations of analysis over surreal number fields}, Mathematics Studies 141, xvi + 373 pp., North Holland, 1987.

\

\noi [ 3 ] {\sc Elwyn R. BERELKAMP, John H. CONWAY, Richard K. GUY}, {\sl Winning ways for your mathematical plays}, vol. 1, Second edition, xix, p. 1-276, A K Peters, Wellesley, Massachusetts, 2001.

\

\noi [ 4 ] {\sc Elwyn R. BERELKAMP, John H. CONWAY, Richard K. GUY}, {\sl Winning ways for your mathematical plays}, vol. 2, Second edition, xviii, p. 277-473, A K Peters, Wellesley, Massachusetts, 2003.

\

\noi [ 5 ] {\sc Elwyn R. BERELKAMP, John H. CONWAY, Richard K. GUY}, {\sl Winning ways for your mathematical plays}, vol. 3, Second edition, xxi, p. 461-801, A K Peters, Wellesley, Massachusetts, 2003.

\

\noi [ 6 ] {\sc Elwyn R. BERELKAMP, John H. CONWAY, Richard K. GUY}, {\sl Winning ways for your mathematical plays}, vol. 4, Second edition, xvi, p. 801-1004, A K Peters, Wellesley, Massachusetts, 2004.

\

\noi [ 7 ] {\sc J. H. CONWAY}, {On numbers and games}, ix + 230 pp.,
Academic Press Inc. 1976, reprinted 1979. 

\

\noi [ 8 ] {\sc J. H. CONWAY}, {\sc All games bright and beautiful}, Amer. Math. Monthly,  {\bf 84}, \no 6 (1977) 417-434. 

\

\noi [ 9 ] {\sc Philip EHRLICH}, {\sl The absolute arithmetic continuum and the unification of all numbers great and small}, The Bulletin of Symbolic Logic, {\bf 18} \no 1 (2012) 1-45.

\

\noi [ 10 ] {\sc Harry GONSHOR}, {\sl An introduction to the theory of surreal numbers}, London Mathematical Society Note Series 110, 192 pp., Cambridge University Press, digitally printed version 2008.

\

\noi[ 11 ] {D.E. KNUTH}, {\sl Les nombres surr\Žels}, {\rm ou comment deux anciens \Žtudiants d\Žcouvrirent les math\Žmatiques pures et v\Žcurent heureux. Une romance math\Žmatique de D. E. Knuth}, Traduction : Daniel E. Loeb et H\Žl\ne Loeb, Original 1974 - Addison Wesley Publishing Company, Traduction March 2, 1997- Loeb, 77 pp.

\

\noi [ 12 ] {Jos\Ž M. PACHECO}, {\sl The Spanish mathematician Norberto Cuesta recovered from oblivion}, Preprint, 11 December 2014.

\

\noi [ 13 ] {\sc Jos\Ž M. PACHECO}, {\sl The life and mathematics of Norberto Cuesta (1907-1989)} 14 October 2016.

\

 \noi [ 14 ] {\sc Simon RUBINSTEIN-SALZEDO,  Ashvin SWAMINATHAN}, {\sl Analysis on surreal numbers}, {\tt arXiv:1307.7392v3}, 19 May 2015.

\

\noi [ 15 ] {\sc Dierk SCHLEICHER, Michael STOLL}, {\sl An introduction te Conway's games and numbers}, {\tt arXiv:0410026v2}, 30 Sep 2005.

\

\noi[ 16 ] {\sc Claus T\O NDERING}, {\sl Surreal numbers - An introduction}, Version 1.7, 31 January 2019.

\

\noi[ 17 ] {\sl Nombre surr\Žel}, WikipediA,  16 avril 2020 \ˆ 11\!:40.

https://fr.wikipedia.org/w/index.php?title

=Nombre surr\Žel\&oldid=169620524.

\

\noi [ 18 ] {\sl Surreal number}, WikipediA, 22 April 2020, , at 08\!:43 (UTC).

https://en.wikipedia.org/w/index.php?title

=Surreal number\&oldid=952449475

\

\noi [ 19 ] BookReview in Bull. A. M. S., {\bf 84} \no 6 (1978) 1328-1336,{\sc Aviezri S. FRAENKEL}, on {\sc Conway} and {\sc Knuth}.

\

\noi [ 20 ] MR886475, M\'arki/Alling.

\

\noi [ 21 ] MR101844, on {\sc Cuesta Dutari}.

\

\noi [ 22 ] MR103838,103839, on {\sc Cuesta Dutari}.

\

\noi [ 23 ] zbMATH 0086.04301, on Cuesta Dutari.

\

\noi [ 24 ] zbMATH 0621.12001, on Alling.

\

\noi [ 25 ] {\sl Sur un article de 1954 sign\Ž N. Cuesta}  :
une traduction, 

\noi {\tt arXiv:2101.05805v1}, 15 Jan 2021.

\

\head{Une derni\re digression}

\

De nos jours, m\me si cela n'est pas tr\s fr\Žquent, il n'est pas rare de voir \Žclore une m\me id\Že math\Žmatique, en plusieurs pays diff\Žrents, quasi simultan\Žment.   Le plagiat, tr\s vilaine chose,  une fois exclus, je pense que cela s'explique simplement de la mani\re suivante, avec la transmission orale !
Plusieurs personnes assistent, en m\me temps, en un m\me lieu, \ˆ un expos\Ž remarquable. Ils sont  impressionn\Žs et, une fois rentr\Žs chez eux, ils leur vient la m\me id\Že ou presque ! Il n'en faut pas davantage. Ce ph\Žnom\ne, en synchronie, existe aussi dans la diachronie, avec la transmission \Žcrite, bien entendu, plus rarement sans doute !

\

\noi Enfin, et pour terminer, je te signale que notre coll\gue, et ami, de Lyon, {\sc Jean-Claude Carrega} que je dois remercier bien vivement ici, ayant lu le brouillon de cette lettre m'a donn\Ž un conseil avis\Ž et m'a encourag\Ž \ˆ la faire publier.

\

\noi Bien \ˆ toi.

\noi Labib

\

\head{Annexes} 

\

\noi Dans ces Annexes,  deux items.

\noi Les pages 4 \ˆ 6 du livre de Conway.

\noi La reproduction de la couverture du num\Žro de la Revista qui contient l'article de Cuesta.

\includegraphics[width=5.3in]{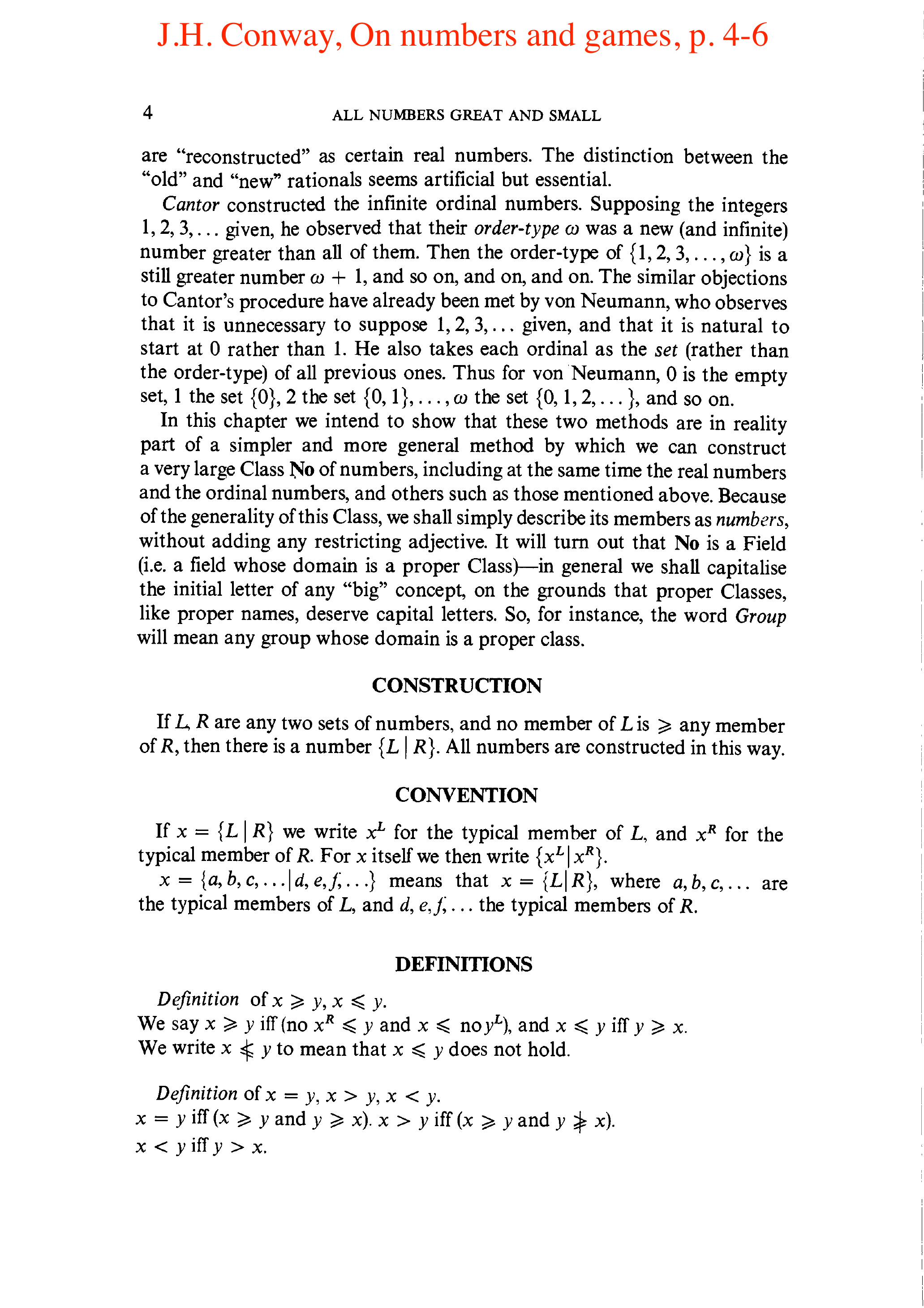}

\includegraphics[width=5.5in]{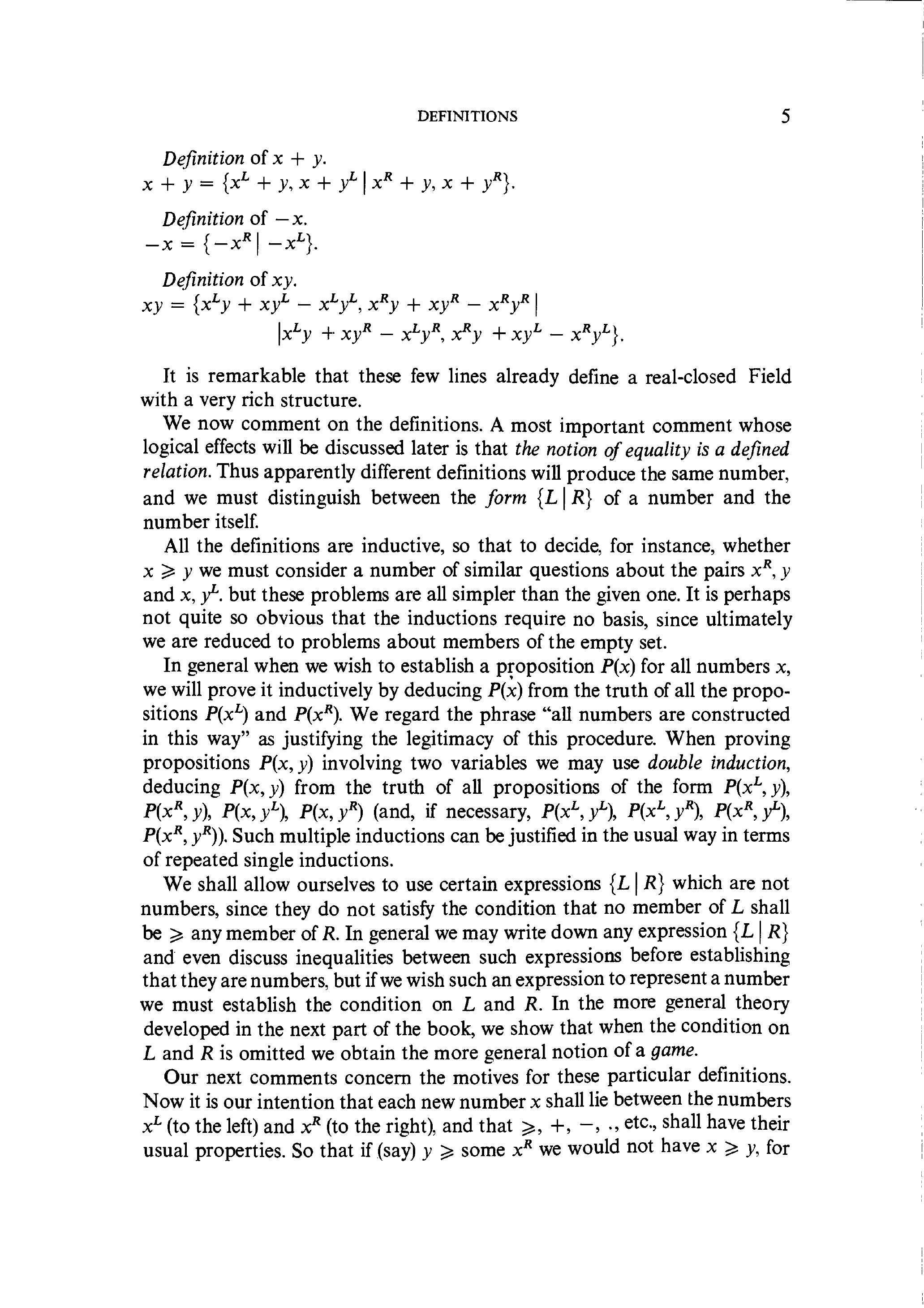}

\includegraphics[width=5.5in]{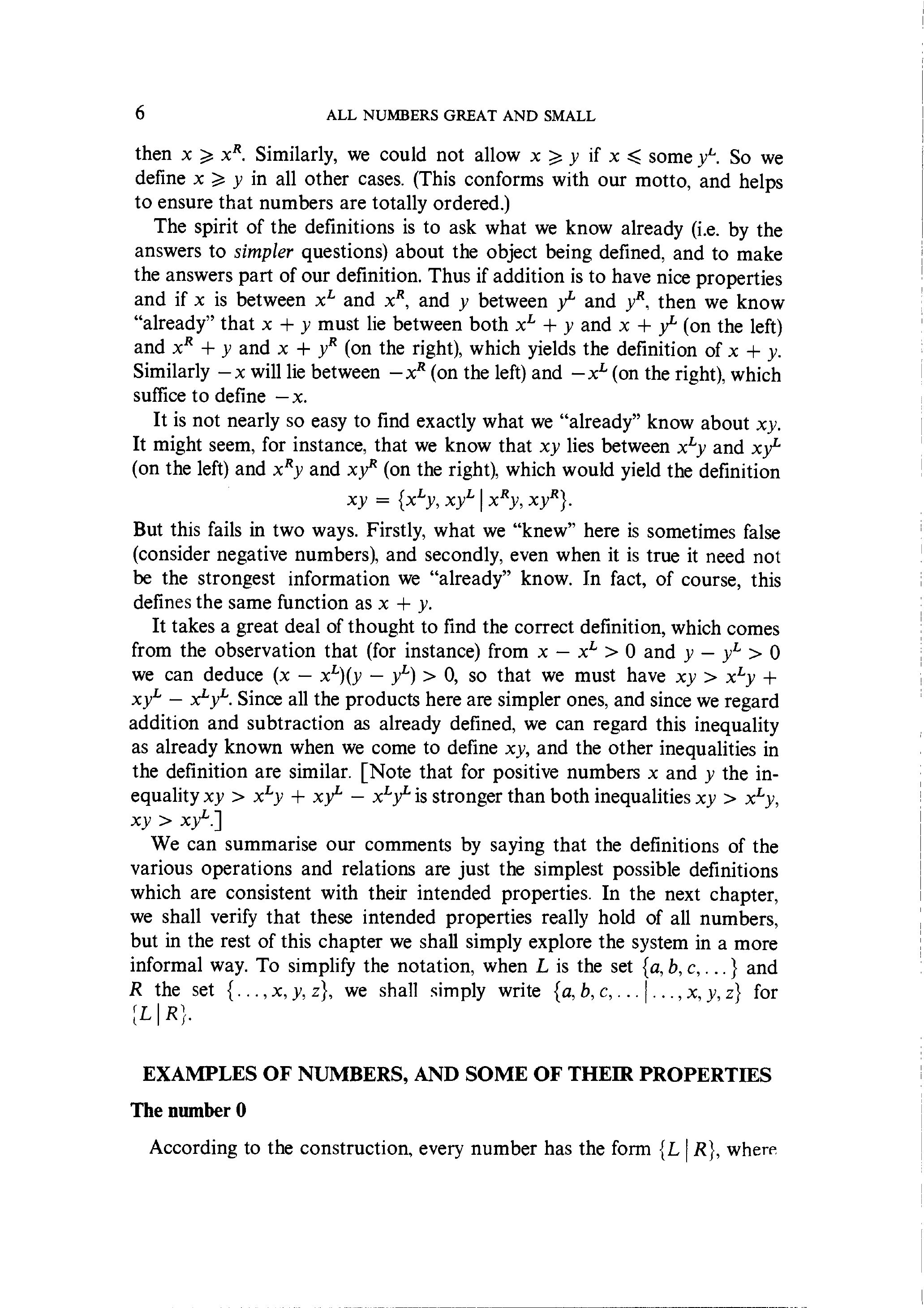}

\includegraphics[width=5.5in]{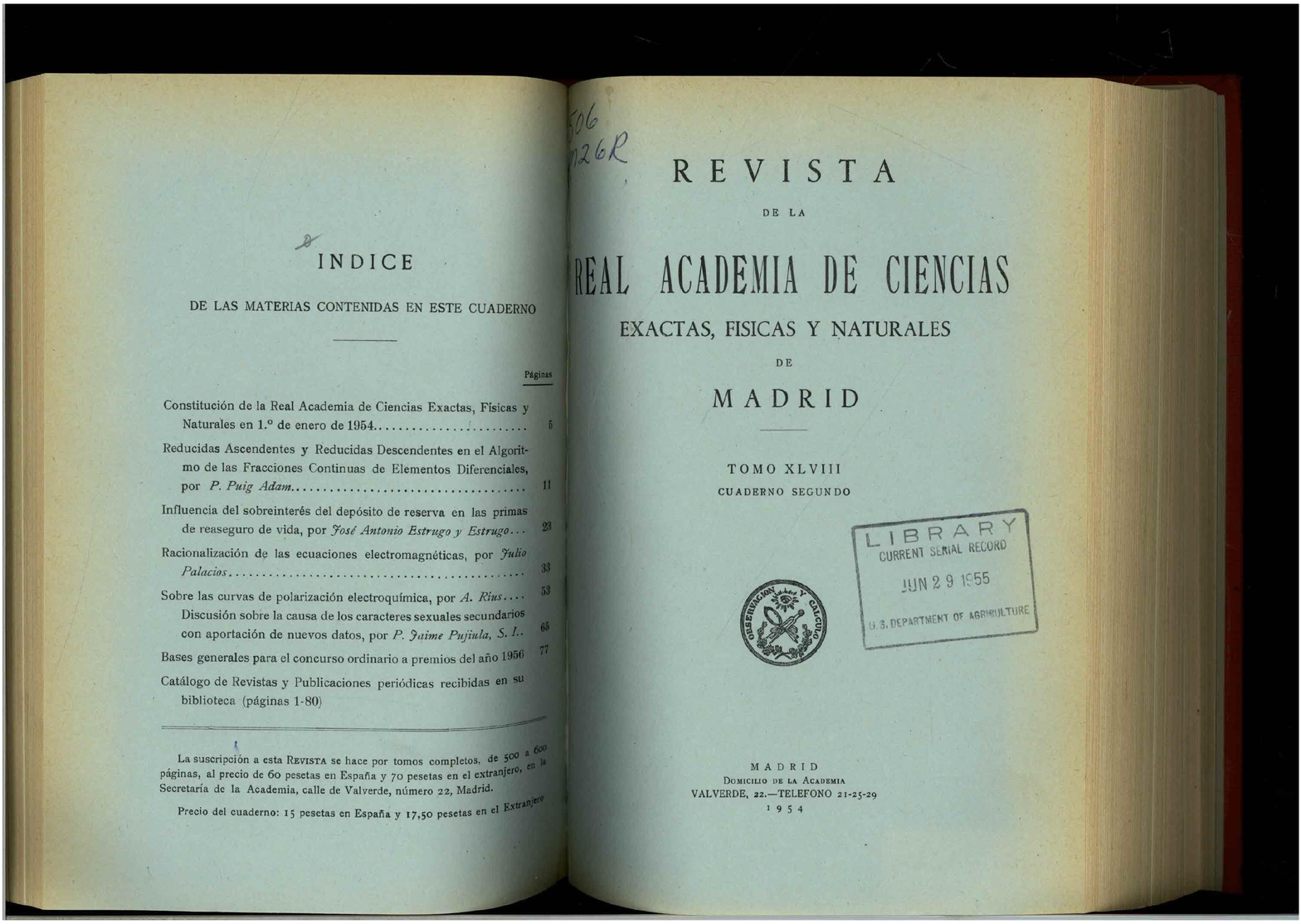}

\

\enddocument

\ 

\

\lopar Le cas o\  $E =\Q$ est l'ensemble des nombres rationnels. Soit $\cal D(\Q)$ l'ensemble des coupures $(A | B)\in \cal C(\Q)$ pour lesquelles il n'y a pas d'\Žl\Žment plus petit que tous les autres dans la partie $B$.  On munit la partie $\cal D(\Q)$ de l'ordre total induit par celui de $\cal C(\Q)$. Obtient ainsi la droite r\Želle $\R$ : si, pour un $x = (A | B) \in \cal D(\Q)$, il y a dans $A$ un \Žl\Žment $q$ plus grand que tous les autres, on identifie $x$ \ˆ $q$. Sinon, $x$ d\Žsigne l'unique nombre r\Žel tel que $A < x < B$.\}

\

\

\

\head{3 \ Une digression} 

\

\

\

\enddocument